\documentclass[letterpaper, 10 pt, draftcls, onecolumn]{ieeeconf} 
\IEEEoverridecommandlockouts 
\usepackage{graphicx}
\usepackage{amsmath}
\usepackage{amssymb}
\usepackage{color}
\usepackage{mathrsfs}
\usepackage{subfigure}

\usepackage{mathtools}
\mathtoolsset{showonlyrefs}

\newtheorem{definition}{Definition}

\newtheorem{theorem}{Theorem}
\newtheorem{remark}{Remark}

\newcommand{\R}{\mathbb{R}}
\newcommand{\N}{\mathbb{N}}
\newcommand{\0}{\emptyset}

\newcommand{\carrier}{\R^n}
\newcommand{\lap}{\textnormal{Lap}}
\newcommand{\diam}{\textnormal{diam}}

\title{\LARGE \bf Differential Privacy for Sets in Euclidean Space}

\author{Matthew T. Hale$^{\star}$%
\thanks{$^{\star}$Department of  Mechanical and Aerospace Engineering, University of Florida, Gainesville, FL, USA.
Email: \texttt{matthewhale@ufl.edu}}%
}

\begin{document}
\maketitle

\begin{abstract}
As multi-agent systems become more numerous and more data-driven,
novel forms of privacy are needed in order to protect data types
that are not accounted for by existing privacy frameworks. In this paper, we
present a new form of privacy for set-valued data which extends
the notion of differential privacy to sets which users want to
protect. 
While differential privacy is typically defined in terms of probability
distributions, we show that it is more natural here to define
privacy for sets over their \emph{capacity functionals},
which capture the probability of a random set intersecting some
other set.
In terms of sets' capacity functionals, we provide a novel 
definition of differential privacy for set-valued data. 
Based on this
definition, we introduce the Laplacian Perturbation Mechanism
(so named because it applies random perturbations to sets),
and show that it provides $\epsilon$-differential privacy 
as prescribed by our definition. 
%While the exponential
%mechanism is often used for privacy of non-numerical data,
%we elaborate upon why its computational complexity makes
%it a less desirable choice here. 
These theoretical results are supported by numerical results,
demonstrating the practical applicability of the developments made. 
\end{abstract}
\section{Introduction}
Multi-agent systems are studied in diverse applications, including 
robotics~\cite{soltero13}, communications~\cite{kelly98,chiang07}, and the 
smart power grid~\cite{caron10,vytelingum10}. 
Agents in these systems work together by exchanging information, 
though in some cases the data that must be shared can be
quite sensitive. For example, it has been shown that the
power usage data shared in smart power grids can reveal a great
deal about users' personal lives, including when a user is at
home or away, awake or asleep, and even what a user is doing inside
his or her home~\cite{edps12,doe10}.
In this and other data-intensive applications, there is therefore a need 
for multi-agent coordination while protecting the data of individuals.  
In response to such privacy concerns, 
privacy tools have been developed that allow for sharing data while
providing privacy guarantees to data-producing entities.

One privacy framework that has seen wide use is differential privacy.
Differential privacy originates in the database literature~\cite{dwork06a,dwork06b,dwork06c}, and its original
formulation keeps individual database entries private whenever a database
is queried. Differential privacy is an attractive means to protect
sensitive data because it provides several strong privacy guarantees.
In particular, it is resilient to post-processing, meaning that post-hoc
analyses do not threaten the privacy protections afforded to sensitive data~\cite{dwork13}, and it is
also robust to arbitrary side information, meaning that learning more about a
data-producing agent does not fully defeat that agent's privacy~\cite{kasiv08}. 
To provide these strong protections to other data types, 
differential privacy has been extended to other settings, including
trajectory-valued data
in~\cite{leny14}, functional data in~\cite{huang15}, and arbitrary metric spaces in~\cite{chatzikokolakis13}. 

Along these lines, in this paper we extend differential privacy to set-valued data.
To consider a concrete example of sensitive set-valued data, we
consider forecasting in financial markets. Given the current market
prices of some collection of related assets, a financial firm can
make estimates of the future prices of these assets, though
each predicted price is often expressed as an interval in which the future
price will lie, rather than as a single point. As a result, predictions about this
collection of assets most naturally take the form of a set~\cite{kaval06}. 
An analyst or government
agency monitoring market behavior may wish to access these predictions
from several firms, though revealing them may be undesirable. Indeed,
revealing these set-valued data can not only make a firm's future intentions
predictable or known, but can also potentially
reveal the manner in which these predictions are made. As a result, 
these set-valued data need to be kept private. 

Beyond this example, set-valued data are also used in a range of diverse applications, including
state estimation~\cite{shamma97,shamma99,ra04}, optimization~\cite{chen06,khan16}, and economics~\cite{kaval06},
and similar privacy concerns exist across these domains.
From these applications, there arises a need to share set-valued data in a manner which makes it private,
while still allowing for meaningful analyses to be performed.
Targeting all of these applications, we develop differential privacy for set-valued data. 

\begin{figure}
\centering     %%% not \center
\subfigure[A connected region in the plane]{\label{fig:a}\includegraphics[width=1in]{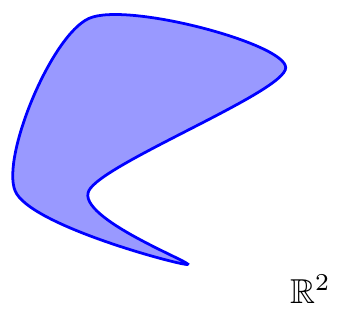}}
\hspace{1pt}
\vrule
\hspace{1pt}
\subfigure[A collection of four isolated points]{\label{fig:b}\includegraphics[width=1in]{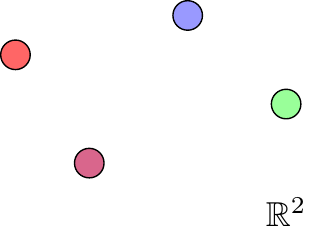}}
\hspace{1pt}
\vrule
\hspace{1pt}
\subfigure[A polytope]{\label{fig:c}\includegraphics[width=1in]{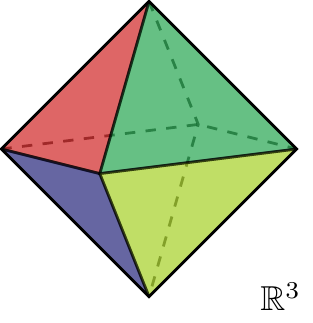}}
\caption{Three examples of sets which can be kept private using the privacy framework we develop. Each set is labeled
with its ambient space in the bottom right corner. We emphasize that these are only
examples, and any closed set in $\R^n$ can be kept private for any $n \in \N$. 
}
\label{fig:sets}
\end{figure}

Several existing works have proposed techniques for keeping sets private, though
many of these works target different kinds of sets than we consider here. 
In particular, emphasis has been placed on privacy for transactional data, which are non-numerical data
that can represent purchases made, and there has been focus on anonymization techniques
for such datasets~\cite{he09,terrovitis08,xu08}, including time-varying datasets~\cite{zhang13},
and 
scalable privacy for high-dimensional sparse datasets of this kind~\cite{ghinita08}. 
In addition to existing work on anonymization,
the work in~\cite{chen11} develops differential privacy for counting queries over transactional
datasets. Motivated by the aforementioned applications, our emphasis in this paper is on 
privacy for subsets of $\R^n$, and we differ from the existing literature precisely because we 
develop differential privacy for subsets of Euclidean space. In general, such sets are not readily
protected by existing methods for set-valued privacy. 
To the best of our knowledge, this is the first work to develop differential privacy for
sets in $\R^n$. 

There are unique challenges in defining privacy for set-valued data, and there are several
desirable characteristics a general theory of private sets should have. 
First, although the framework we develop should certainly accommodate sets with non-empty interior,
it should also apply just as well to singletons or sets comprised by isolated points, simply because
doing so significantly broadens the scope of this work. 
Second, the shape of the set should be preserved. In the above financial forecasting example,
the exact set should be kept private, though if it is a cube, parallelipiped, sphere, or other
shape, that aspect of the set should be preserved in order to keep the privatized data in a 
state amenable to analysis. Third, the privacy framework developed should be
computationally simple so that it can be seamlessly integrated into a range of applications
with minimal overhead and without the need for specialized computing hardware. 

To develop differential privacy for set-valued data, we will draw from existing work
on the theory random sets. A general theory of random sets has been developed over time, driven
by a range of applications from geology~\cite{matheron70} to microscopy~\cite{matheron72}. 
The work in \cite{molchanov06} appears to provide the most
comprehensive reference on the subject, and we therefore use it as our starting point
for random sets. 
The notion of a random set is distinct from that
of a fuzzy set because a fuzzy set allows for degrees of inclusion in a set, 
while a random set is simply a random variable whose outcomes are ordinary sets, and we
favor this framework because privatized sets should still be ordinary sets in 
order to be useful in applications dependent upon set-valued data. 
The work in \cite{molchanov06} develops a theory of random closed sets, 
and we build upon this work to provide privacy for sets in a manner
that meets the three criteria described above. 
%and this immediately accommodates
%singleton sets, the first requirement stated above, and it will be shown 
%in Section~\ref{sec:mech} that building upon this framework allows
%for preservation of the shape of a set and computational simplicity,
%which were the other two requirements stated above.

Beyond our own requirements,
there are several challenges unique to privacy of set-valued data, including 
substantial difficulties in even defining probability distributions over 
rich spaces of sets. 
We address these difficulties in our definition of privacy for sets, and
one contribution of this paper is the privacy definition itself.
Our definition deviates from the usual definition
of differential privacy because it is not in terms of probability distributions. Instead,
our definition is in terms of the \emph{capacity functionals} of random sets, 
which give the probability of a random set intersecting another set, precluding the
need for using explicit probability distributions. 
This definition in turn lets
us consider the events of interest in a significantly simpler setting; it
will be shown that the definition here is more appropriate for set-valued privacy and does not
incur any meaningful losses in privacy or applicability in exchange for its simplicity. 
The definition used here is an adaptation of existing work on differential privacy
on metric spaces, and differential privacy as we define it is a direct generalization
of the work in \cite{chatzikokolakis13} to the case of set-valued data. In terms 
of this definition, we define the Laplacian Perturbation Mechanism for privacy of sets,
so named because it applies a random perturbation to a set in order to make it private. 

The rest of the paper is organized as follows. Section~\ref{sec:background} reviews the
necessary background on the theory of random sets and spaces of sets. 
Then Section~\ref{sec:privacy} provides exposition on differential privacy,
including its generalization to metric spaces. Next, Section~\ref{sec:mech} provides
our definition for set-valued privacy, and defines a privacy mechanism that satisfies
this definition. After that, Section~\ref{sec:numerical} provides numerical results
showing the applicability of the theoretical developments we make. Finally,
Section~\ref{sec:conclusion} concludes the paper.

\section{Background on Random Sets and Spaces of Sets} \label{sec:background}
In this section we provide the necessary background on random sets
and the associated definitions needed to study them. Unless otherwise indicated, the exposition
in this section draws from Chapter~1 of \cite{molchanov06}, and we refer
the reader to that reference for an extensive treatment of this subject. 
We begin by reviewing the necessary background on the theory of random sets, and then we discuss
the metric which will be used in Section~\ref{sec:privacy} to define differential privacy
for sets. 

\subsection{Definition of a Random Set}
As described above, we work exclusively with sets in $\carrier$ because
the applications of interest generate such sets. 
In this setting, we use
$\mathcal{G}$ to denote all open subsets of $\carrier$, $\mathcal{F}$ for
all of its closed subsets, and $\mathcal{K}$ for all of its compact subsets. We also fix a
complete probability space $(\Omega, \mathfrak{F}, \mathbb{P})$, which will be used
throughout when defining random objects. As noted in the introduction, we are interested
in creating a general theory of privacy for sets which incorporates singleton sets (and unions thereof),
and one way to do so is to develop privacy for arbitrary closed sets. Along these lines, 
the work in~\cite{molchanov06} provides a rich theory of random closed sets, and we therefore
use it as our starting point for developing set-valued privacy. 
Toward doing so, we now 
present the definition of a random closed set.

\begin{definition} \label{def:randset} \emph{(Random closed set)}
A random closed set is a map $X : \Omega \to \mathcal{F}$ satisfying
\begin{equation}
\{\omega \in \Omega \mid X(\omega) \cap K \neq \0\} \in \mathfrak{F}
\end{equation}
for all $K \in \mathcal{K}$. \hfill $\triangle$
\end{definition}

Intuitively, Definition~\ref{def:randset} says that $X$ is a random set
if one can observe $X$ and always determine if $X$ does or does not intersect
any compact set $K \in \mathcal{K}$. An equivalent characterization that
highlights key probability-theoretic points is as follows. For all
$K \in \mathcal{K}$, define the families of sets
\begin{align}
\mathcal{F}_K &:= \{F \in \mathcal{F} \mid F \cap K \neq \0\} \\
\mathcal{F}^K &:= \{F \in \mathcal{F} \mid F \cap K = \0\},
\end{align}
and for all $G \in \mathcal{G}$ define
\begin{equation}
\mathcal{F}_G := \{F \in \mathcal{F} \mid F \cap G \neq \0\}. 
\end{equation}

The \emph{Fell topology} on $\mathcal{F}$ is then generated by
the open sets $\mathcal{F}_G$ and $\mathcal{F}^K$ for
all $G \in \mathcal{G}$ and $K \in \mathcal{K}$. 
Using $\mathscr{B}(\mathcal{F})$ to denote the 
Borel $\sigma$-algebra generated by the Fell topology on $\mathcal{F}$,
 we note that
$\mathscr{B}(\mathcal{F})$ is also generated by
the open sets $\mathcal{F}_K$ for all $K \in \mathcal{K}$. 
This characterization of $\mathscr{B}(\mathcal{F})$  leads to the
alternative definition of a random closed set as a map $X: \Omega \to \mathcal{F}$
where
\begin{equation}
X^{-1}(\mathcal{X}) = \{\omega \in \Omega \mid X(\omega) \in \mathcal{X}\} \in \mathfrak{F}
\end{equation}
for all $\mathcal{X} \in \mathscr{B}(\mathcal{F})$. In words, this definition
specifies the precise sense in which a random closed set is measurable with respect to the
Borel $\sigma$-algebra $\mathscr{B}(\mathcal{F})$.

Differential privacy is typically defined with respect to a {$\sigma$-algebra},
and frequently the Borel $\sigma$-algebra is favored in such definitions
because it is rich enough to capture all events of interest. 
While $\mathscr{B}(\mathcal{F})$ is the Borel $\sigma$-algebra over
the Fell topology in this setting, it is 
large enough that it is
difficult to explicitly define a probability measure over it, and, as a result,
probability distributions are often not directly
used in the study of random sets~\cite[Page 10]{molchanov06}. Instead, using the aforementioned
fact that the collection $\{\mathcal{F}_K\}_{K \in \mathcal{K}}$ also
generates $\mathscr{B}(\mathcal{F})$, attention is given specifically to the events
captured by the sets $K \in \mathcal{K}$. In terms of this collection, we
now define the capacity functional of a random set, which 
is significantly simpler than the probability distribution of a random set. 

\begin{definition} \label{def:capfunc} \emph{(Capacity functional)}
The capacity functional of a random set $X$ is denoted ${T_X : \mathcal{K} \to [0, 1]}$ and is defined by
\begin{equation} 
T_X(K) := \mathbb{P}[X \cap K \neq \0]
\end{equation}
for all $K \in \mathcal{K}$. \hfill $\triangle$
\end{definition} 

In words, the capacity functional gives the probability that the random set $X$ intersects
the compact set $K \in \mathcal{K}$. 
As noted above, we choose to work with capacity functionals rather
than probability distributions because doing so is significantly simpler. 
However, this simplicity does not result in any significant loss of
precision in terms of observable events; Definition~\ref{def:capfunc} lets us examine the
probability of $X \in \mathcal{F}_K$ for any $K \in \mathcal{K}$, and, because
these sets in turn generate $\mathscr{B}(\mathcal{F})$, we see that
most events that would be captured by a set's probability distribution over $\mathscr{B}(\mathcal{F})$
 are also captured by its
capacity functional. 

Furthermore, 
Choquet's Theorem \cite[Theorem 1.13]{molchanov06} provides that, under mild conditions, a capacity
functional of the form above uniquely determines a random set, demonstrating that no specificity is lost
in exchange for the simplicity of the capacity functional. 
With this idea in mind, we formally define differential privacy for sets in Section~\ref{sec:privacy}
in terms of the capacity functionals of the sets of interest. Before doing so, we next define
the metric which we use in order to formally
define differential privacy.

\subsection{The Hausdorff Metric}
The introduction discussed the use of metric differential privacy~\cite{chatzikokolakis13} in this setting, and one
naturally must choose which metric to use based upon the needs of privacy in this context.
Differential privacy on metric spaces is designed to make nearby objects in the space
approximately indistinguishable, and the choice of metric provides an opportunity to
define ``nearby'' in the correct way. In the applications of interest, two sets are
typically regarded as nearby if the points they contain are all nearby. Therefore, while
there are other choices for metrics on spaces of sets, e.g., the Fr\'{e}chet-Nikodym
metric \cite[Section 1.12(iii)]{bogachev07} which gives the distance between two sets as the volume of their mutually
non-overlapping regions, we choose to use the Hausdorff metric, as it accurately
captures the notion of closeness required here. We formally define it now in the form
in which we use it below. 

\begin{definition} \label{def:hausdorff} (\emph{\cite{hausdorff62}; Hausdorff Metric})
Let $\|\cdot\|_1$ denote the $1$-norm on $\mathbb{R}^n$. 
For two sets $S_1, S_2 \subset \mathbb{R}^n$, the \emph{Hausdorff distance} between them is defined as
\begin{equation}
d_H(S_1, S_2) := \max\left\{\sup_{s_1 \in S_1} \inf_{s_2 \in S_2} \|s_1 - s_2\|_1, \sup_{s_2 \in S_2} \inf_{s_1 \in S_1} \|s_1 - s_2\|_1\right\},
\end{equation}
where $d_{H}$ is called the \emph{Hausdorff metric}. \hfill $\triangle$
\end{definition}

Below we use $d_H$ to define privacy in a way which modulates the amount of noise added 
as a function of the distance between sets. 
%While either term in the maximum defining $d_H$ has some properties of a metric, taking
%the maximum of these two terms is required for $d_H$ to have the symmetry and homogeneity
%properties required to be a metric. 
We will also use the related notion of the 
diameter of a set, which is given by
\begin{equation}
\diam(S) = \sup_{s_1, s_2 \in S} \|s_1 - s_2\|_1
\end{equation}
for a set $S \subset \R^n$. 
Having defined the necessary aspects of spaces of sets, we proceed to formally define
differential privacy in the next section.

\section{Differential Privacy} \label{sec:privacy}
In this section, we review the necessary background on the
theory of differential privacy. We first briefly review
differential privacy in its usual form, and then elaborate
on the extension of this definition to arbitrary metric spaces.
We then state our definition of differential privacy for
set-valued data.

\subsection{Background on Differential Privacy}
Differential privacy was originally defined for database queries in \cite{dwork06a,dwork06b,dwork06c}, and 
its underlying principle is simple: for a database query $q$, noise is added so that
nearby databases produce query responses that are statistically close. More concretely,
let $\mathsf{D}$ be the set of databases of interest (containing ages, salaries, etc.).
For any choice of privacy parameter $\epsilon \geq 0$, privacy is enforced
by a \emph{mechanism}, which is a randomized function meant to approximate exact queries
by privatized queries, with the level of privacy specified by $\epsilon$. 
Formally, a mechanism $M$ on $\mathsf{D}$ is $\epsilon$-differentially private if, for all databases
$D_1, D_2 \in \mathsf{D}$ differing by a single entry, we have
\begin{equation}
\mathbb{P}[M(D_1) \in A] \leq e^{\epsilon}\mathbb{P}[M(D_2) \in A]
\end{equation}
for all measurable sets $A$ in the range of $M$. The purpose of protecting single
changes in database values is to protect individuals' data while still providing
accurate responses to queries of groups of users. 
The privacy parameter $\epsilon$ determines how private data must be kept, with smaller
values of $\epsilon$ giving stronger privacy protections, and typical values
range from $0.1$ to $\log 3$. 
Given $\epsilon$, the above definition
specifies the exact relationship required of private queries of $D_1$ and $D_2$,
and the strong privacy protections provided by this definition have been well-studied~\cite{dwork13}. 

Indeed, several key properties of differential privacy can be derived from this
definition. First, differential privacy is resilient to post-processing, so that
$M$ being differentially private implies that $f \circ M$ is as well for
any deterministic $f$ such that the composition is well-defined. This allows
differentially private data to be freely processed without threatening the
privacy guarantees afforded to the underlying sensitive data. 
Second, differential privacy is robust to side information, and
this property implies that learning some database entries does not, 
in general, fully defeat differential privacy \cite{kasiv08}. 
These properties make differential privacy an appealing choice
for protecting data because they allow for the release
of private data without concern for what an eavesdropper will do with
this data or what else they may learn in the future. 

\subsection{Differential Privacy on Metric Spaces}
In this paper, we do not use the above definition for differential privacy, instead opting
for a generalized definition for arbitrary metric spaces. 
The ordinary definition of differential privacy requires one to bound the
\emph{sensitivity} of the query $q$, which is an upper bound on how much
$q(D_1)$ differs from $q(D_2)$ when $D_1$ and $D_2$ differ by a single element.
For set-valued data, queries can themselves be set-valued or scalar-valued, and
the wide range of possible query outputs makes it difficult to compute
general-purpose sensitivity bounds. To avoid doing so, we use the formulation
of differential privacy for metric spaces. 

To define this type of differential privacy, we
consider an arbitrary space $\mathcal{S}$, 
and we use a
metric over this space, defined as $d : \mathcal{S} \times \mathcal{S} \to \R$. 
A mechanism $M : \mathcal{S} \to \mathcal{S}$ is then said to be $\epsilon d$-differentially private if
\begin{equation} \label{eq:metricdp}
\mathbb{P}[M(S_1) \in A] \leq e^{\epsilon d(S_1, S_2)} \mathbb{P}[M(S_2) \in A]
\end{equation}
for all $S_1, S_2 \in \mathcal{S}$ and all measurable sets $A$ in the range of $M$.  

In addition to eliminating the need for sensitivity bounds, this definition has
the feature that it requires some level of indistinguishability between
all choices of $S_1$ and $S_2$, regardless of how different they are. 
In particular, while ordinary differential privacy only pertains to
databases differing by a single entry, this definition applies to all pairs
of input data and scales privacy based on the extent to which these
data differ, providing some level of privacy in all cases. 
With this definition in hand, we now formally define differential privacy
for set-valued data. 

\subsection{Differential Privacy for Set-Valued Data}
We now state the definition for differential privacy that
is used throughout the remainder of the paper. This definition
builds upon the aforementioned definition of privacy in metric spaces, though
it adds an additional term which accounts for this definition being in terms
of capacity functionals rather than probability distributions. We first
state the definition itself and, below that, we discuss the need for
this difference. 

\begin{definition} \label{def:main} \emph{(Differential privacy for sets)}
A mechanism $M : \mathcal{F} \to \mathcal{F}$ is $\epsilon d_H$-differentially private
if, for any two sets $X, Y \in \mathcal{F}$, 
\begin{equation}
T_{M(X)}(K) \leq e^{\epsilon \diam(K)}e^{\epsilon d_H(X, Y)}T_{M(Y)}(K)
\end{equation}
for all $K \in \mathcal{K}$, where $d_H$ is
the Hausdorff metric on $\mathcal{F}$, and $T_{M(X)}$ and $T_{M(Y)}$ are the 
capacity functionals of the sets $M(X)$ and $M(Y)$, respectively. \hfill $\triangle$
\end{definition}

\begin{remark} \label{rem:whydef} 
The presence of the multiplicative term $e^{\epsilon \diam(K)}$ causes this
definition for differential privacy to differ from the framework established
in \cite{chatzikokolakis13}. 
The reason for this change stems from using capacity functionals
instead of probability distributions, 
though its interpretation 
is intuitive: the capacity
functionals of $M(X)$ and $M(Y)$ are parameterized by the choice of set $K \in \mathcal{K}$, and
as $K$ changes, the privacy needs of the underlying data must change in response. In particular,
if $K$ is a small set, then determining that $M(X)$ intersects $K$ but $M(Y)$ does not
allows for fine differences between $M(X)$ and $M(Y)$ to be discerned. 
Thus, in this case, the definition of privacy should be informed by the small size
of $K$, and therefore the term $e^{\epsilon \diam(K)}$ is incorporated, requiring
further closeness between the capacities of $M(X)$ and $M(Y)$ when $K$ is small. Similarly,
if $K$ is a large set, then knowing that both $M(X)$ and $M(Y)$ intersect $K$ does not
threaten privacy much, simply because many sets will intersect a large $K$. Accordingly,
the term $e^{\epsilon \diam(K)}$ encodes this fact by requiring weaker privacy guarantees
for very large choices of $K$. \hfill $\lozenge$
\end{remark}

\begin{remark} \label{rem:asymptotically}
Asymptotically, as $K$ becomes very large, the term $e^{\epsilon \diam(K)}$ goes to infinity,
encoding the fact that no privacy protections would be needed for an infinite set.
Indeed, so many sets
would intersect such a $K$ that there would be no privacy loss in detecting its
intersection with another set, and this intuition is captured in the $e^{\epsilon \diam(K)}$ term. 
Similarly, if $K = \{\xi\}$ is a singleton, then the capacity
functional of a set determines the probability that a set contains $\xi$. Here, there are significant
privacy concerns because the knowledge that $\xi \in M(X)$ and $\xi \not\in M(Y)$ allows for
very fine distinctions to be made between $M(X)$ and $M(Y)$, 
leading to substantial losses of privacy. 
In such a case, $\diam(K) = 0$ and $e^{\epsilon \diam(K)} = 1$, reducing the above privacy definition
to the one seen for metric spaces in Equation~\eqref{eq:metricdp}, and enforcing stronger privacy
requirements. 
In light of this fact, the privacy definition
given in Definition~\ref{def:main} can be seen to be a direct generalization of the one
in \cite{chatzikokolakis13}. \hfill $\lozenge$
\end{remark}

Having defined privacy for set-valued data, the next section defines the Laplacian Perturbation
Mechanism and proves that it provides $\epsilon d_H$ privacy as defined by Definition~\ref{def:main}. 

%\red{Talk about why $\diam(K)$ is scaled by $\epsilon$? Talk about why it's exponentiated?} 

\section{Enforcing Differential Privacy for Set-Valued Data} \label{sec:mech}
In this section we define the 
Laplacian Perturbation Mechanism and show that it provides
differential privacy as defined in Section~\ref{sec:privacy}. Then we discuss
why we do not use the exponential mechanism, despite its common use for
non-numerical data. 
Below, we use the notation
\begin{equation}
X + a := \{x + a \mid x \in X\}
\end{equation}
for all $X \subseteq \R^n$ and all $a \in \R^n$. 

\subsection{The Laplacian Perturbation Mechanism}
In this subsection we define the Laplacian Perturbation Mechanism for 
privacy of set-valued data, which, as its name suggests, applies
perturbations to sets which are randomly drawn from a Laplace distribution. 
The scalar Laplace distribution has
probability density
\begin{equation}
\lap(\mu, b) := \frac{1}{2b} \exp\left(-\frac{|x-\mu|}{b}\right),
\end{equation}
where $\mu \in \R$ is the mean of the distribution and $b \in \R$ is its scale parameter. 
Using this distribution we have the following result. 

\begin{theorem} \label{thm:lapmain} 
The Laplacian Perturbation Mechanism in $\carrier$ defined by
\begin{equation}
M(X) := X + w, \,\, \textnormal{ where } w \sim \lap\left(0, \frac{1}{\epsilon}\right)^n
\end{equation}
is $\epsilon d_H$-differentially private for all sets $X \in \mathcal{F}$. 
\end{theorem}
\emph{Proof:} 
Consider an arbitrary compact set $K \in \mathcal{K}$, along with
two sets $X, Y \in \mathcal{F}$. We find that, for a vector
$w \in \R^n$ with $w_i \sim \lap(0, 1/\epsilon)$, expanding definitions gives
\begin{align}
\frac{T_{M(X)}(K)}{T_{M(Y)}(K)} &= \frac{T_{X+w}(K)}{T_{Y+w}(K)} \\ 
                                &= \frac{\mathbb{P}[X+w \cap K \neq \0]}{\mathbb{P}[Y+w \cap K \neq \0]} \\
                                &= \frac{\mathbb{P}[w = x - k \textnormal{ for some } x \in X, \, k \in K]}{\mathbb{P}[w = y - \xi \textnormal{ for some } y \in Y, \, \xi \in K]}.
\end{align}
We emphasize that $k$ and $\xi$ need not have
any relationship because the capacity functionals here capture the probability of intersecting with the set $K$, though this analysis does not require
that $M(X)$ and $M(Y)$ intersect $K$ in the same place, allowing $k$ and $\xi$ to differ. Continuing, we examine these probabilities component-wise, and use
the definition of the Laplace distribution given above, to find
\begin{align}
\frac{\mathbb{P}[w=x-k\textnormal{ for some }x\in X, \, k\in K]}{\mathbb{P}[w=y-\xi\textnormal{ for some }y\in Y, \, \xi\in K]} &= \prod_{i=1}^{n} \frac{\mathbb{P}[w_i = x_i - k_i, \, x \in X, \, k \in K]}{\mathbb{P}[w_i = y_i - \xi_i, \, y \in Y, \, \xi \in K]} \\
                &= \prod_{i=1}^{n}\frac{\exp(-\epsilon|x_i-k_i|)}{\exp(-\epsilon|y_i-\xi_i|)}, \textnormal{ for } x \in X, \, y \in Y, \, k, \xi \in K \\
                &= \exp\left(-\epsilon\sum_{i=1}^{n}|x_i - k_i| + \epsilon\sum_{i=1}^{n}|y_i - \xi_i|\right), \, x \in X, \, y \in Y, \, k, \xi \in K \\
                &=\exp\Big(\epsilon(\|y - \xi\|_1 - \|x - k\|_1)\Big), \textnormal{ for } x \in X, \, y \in Y, \, k, \xi \in K \\
                &\leq\exp\Big(\epsilon \|(y - x) + (k - \xi)\|_1\Big), \textnormal{ for } x \in X, \, y \in Y, \, k, \xi \in K,
\end{align}
which follows from the reverse triangle inequality. Continuing, we use the ordinary triangle inequality to find
\begin{equation}
\exp(\epsilon \|(y - x) + (k - \xi)\|_1) \leq \exp(\epsilon\|k - \xi\|_1)\exp(\epsilon\|x - y\|_1) \textnormal{ for } x \in X, \, y \in Y, \, k, \xi \in K,
\end{equation}
whereupon we use ${\|k - \xi\|_1 \leq \diam(K)}$ and $\|x - y\|_1 \leq d_H(X, Y)$ to conclude. 
\hfill $\blacksquare$

Theorem~\ref{thm:lapmain} shows that a set can be kept differentially private by applying a random perturbation
in each coordinate. The introduction specified that the privacy mechanism developed in this paper should be computationally
simple, and we see that this is indeed the case as keeping a set in $\R^n$ private can be achieved by generating
$n$ random numbers. Next we comment on why we deliberately do not use the exponential mechanism for privacy
in this setting. 

\subsection{On Avoiding the Exponential Mechanism}
Perhaps the best known privacy mechanism for non-numerical data is the exponential
mechanism. In using it, output data are assigned a quality score, which encodes how well
each possible output approximates the true output of a query. Then,
private outputs of non-numerical queries are generated in a way that outputs
high-utility outputs with high probability and low-utility outputs with
low probability. 
To define the exponential mechanism 
for closed sets in $\R^n$, one would need to define a utility function
\begin{equation}
u : \mathcal{F} \times \mathcal{F} \to \R
\end{equation}
which provides a quantitative measure of the utility of outputting the set
$\bar{S} \in \mathcal{F}$ when the input set is $S \in \mathcal{F}$. 

A known theoretical hurdle \cite[Page 38]{dwork13} 
in implementing the exponential mechanism is
the complexity associated with generating quality scores for all
possible input/output pairs, and this can be prohibitive for some data types.
Here, this hurdle precludes the use of the exponential mechanism
for the same reason we are precluded from working directly with probability
distributions: the space of all subsets is simply
too large to assign each subset a quality score, and
thus we avoid doing so and likewise avoid using the exponential mechanism. 
However, as was shown above, we have
successfully used the Laplace mechanism for privacy, meaning
that differential privacy can be provided for sets despite the inability to use
the exponential mechanism. 

Next, we give numerical implementations of the Laplacian Perturbation Mechanism
to show its use in practice.

%\red{Relationship to random intervals?} 
%\red{At some point, let's have a section about keeping intervals private
%using ordinary differential privacy, and then let's show that that basically
%collapses intervals to sets with different measure. Therefore, something is needed
%which respects the fact that the data are sets and preserves this nature while 
%still allowing for privacy.} 

\section{Numerical Results} \label{sec:numerical}
In this section we present two sets of numerical results
implementing the Laplacian Perturbation Mechanism
defined in Theorem~\ref{thm:lapmain}. In order to generate
results that are easily plotted, we do two numerical
runs, one in $\R^2$ and one in $\R^3$. 

\subsection{Numerical Results in $\R^2$} \label{ss:results1}
To demonstrate the effects of the Laplacian Perturbation Mechanism
in $\R^2$, we implement it on a set of four points. In particular, 
we begin with the set
\begin{equation}
X = \{(0, 0), (0, 1), (1, 0), (1, 1)\} \subset \R^2
\end{equation}
and generate four outputs of the Laplacian Perturbation Mechanism
with this set as the input. Here the privacy parameter was
selected to be $\epsilon = 1$. In this case, each of the four
numerical runs implemented privacy by adding noise $w$, where
\begin{equation}
w \sim \lap(0, 1)^2.
\end{equation}

\begin{figure}
\centering
\includegraphics[width=3.3in]{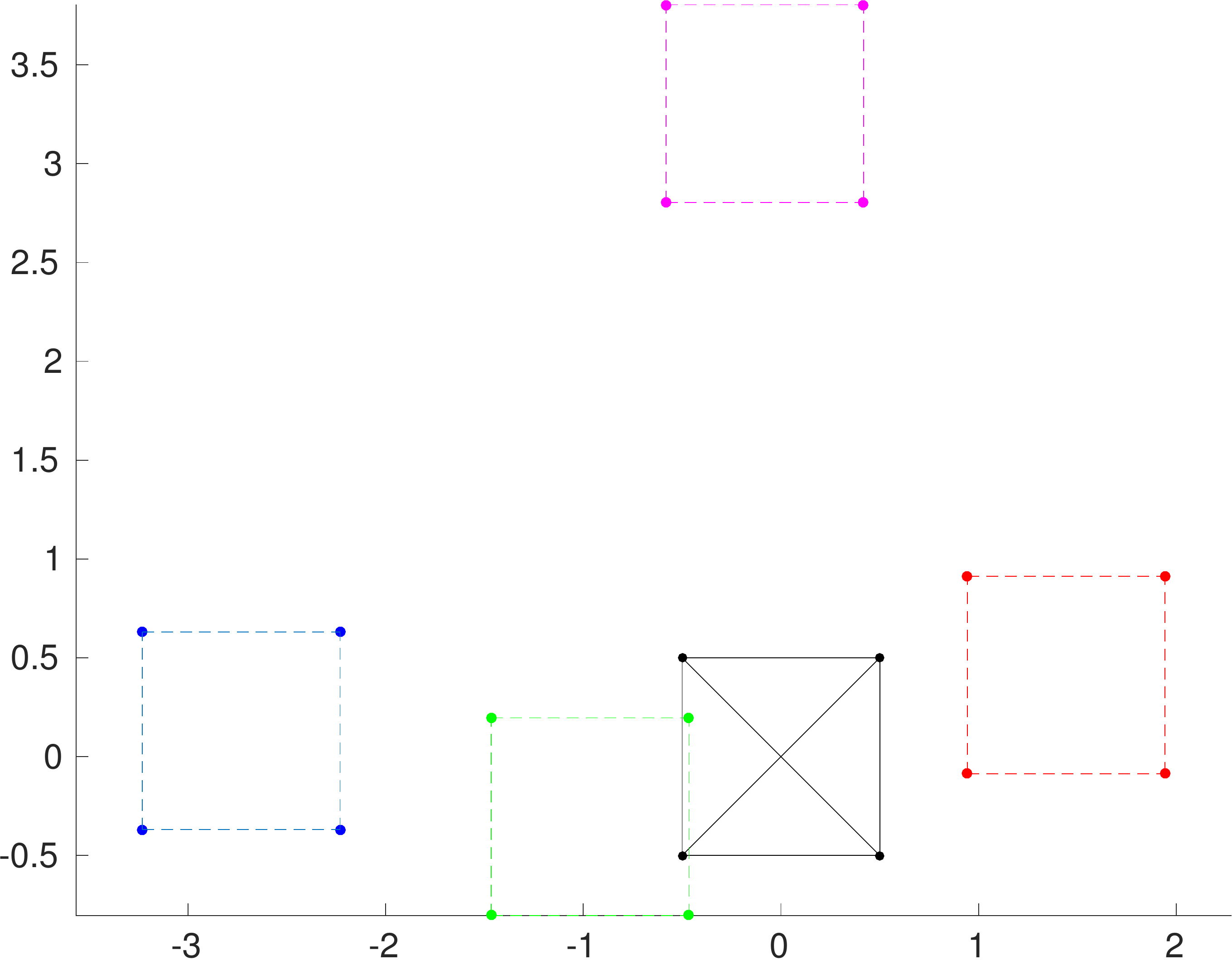}
\caption{The numerical results of four runs of the Laplacian Perturbation
Mechanism in $\R^2$. The original set to keep private consists of four
isolated points, and these points are connected by a square with an
``X'' in the middle to indicate that it is the sensitive data we want
to protect. Each query of this set is shown with its four points connected
by dotted lines to indicate which points belong to the same set. Here we can
visually verify that privacy is provided while requiring only mild
perturbations to the original set, indicating that private set-valued
data will still be useful in the computations required by
the applications of interest. 
}
\label{fig:points}
\end{figure}

The results from these runs are shown in Figure~\ref{fig:points}. The original
dataset is shown in the black square, with the points connected by solid
lines and an $X$-shape in the center. We emphasize that each set only consists
of the vertices of each square, and we connect points in this way
simply to make clear which points
belong to the same set. We see that the private outputs
of the Laplacian Perturbation Mechanism are able to preserve privacy
of the set $X$ while still allowing its outputs to remain close to $X$
itself, indicating that privacy can be implemented while still retaining
reasonably accurate set-valued data. 

\subsection{Numerical Results in $\R^3$}
\begin{figure}
\centering
\includegraphics[width=3.3in]{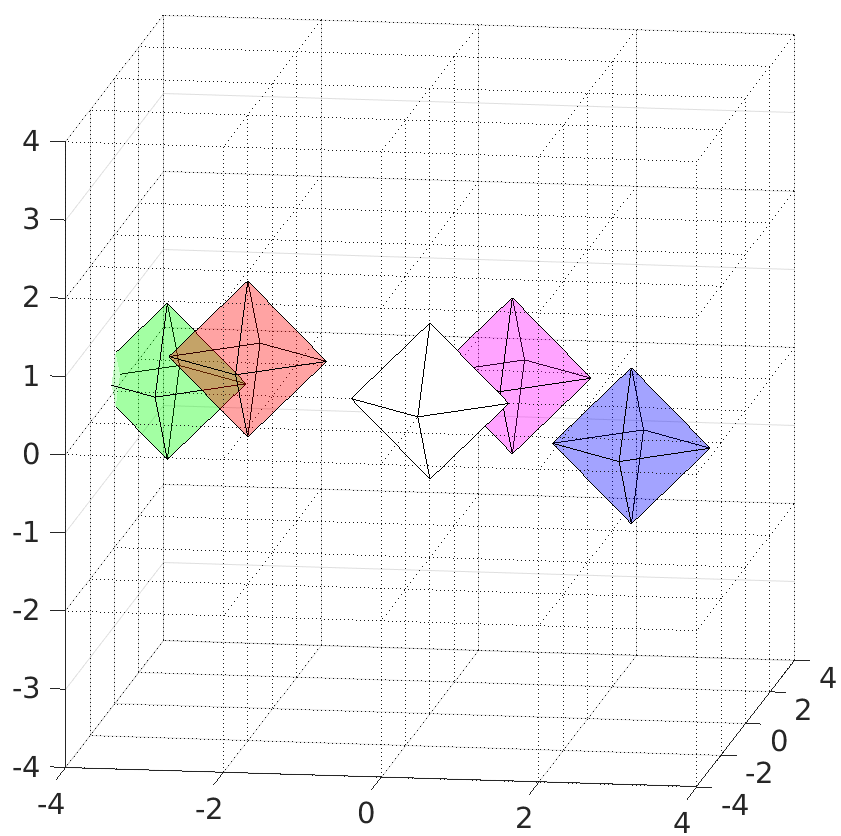}
\caption{The numerical results of four applications of the Laplacian
Perturbation Mechanism in $\R^3$. The original sensitive data is the
polytope in the center, and the outputs of the four numerical
runs are the other four polytopes. As in Figure~\ref{fig:points}, we see
that privacy is attained with reasonably sized perturbations, allowing
for accurate responses to private queries of sensitive set-valued data. 
}
\label{fig:poly}
\end{figure}

To further illustrate the use of the Laplacian Perturbation Mechanism, we 
give numerical results for a set contained in $\R^3$. The set we begin
with is the polytope in the center of Figure~\ref{fig:poly},
shown in white. As above, four total numerical runs were conducted,
each with $\epsilon = 1$, giving
\begin{equation}
w \sim \lap(0, 1)^3
\end{equation}
in each case. 

The results of these runs are shown in Figure~\ref{fig:poly}, where
the center polytope is the original set which we intend
to keep private. Each of the other four polytopes results from
a numerical run we conduct, and it can be readily verified
in Figure~\ref{fig:poly} that the random displacements needed
for differential privacy are on the order of the diameter
of the set itself, indicating that set-valued data can be
protected with reasonable random perturbations. 

The results of this subsection
and those above in Section~\ref{ss:results1} indicate that the privacy
framework we develop meets the three criteria given in
the introduction, namely, that collections of singletons
should be accommodated, that the shape of the set should
be preserved, and that the framework should be computationally
simple. Moreover, we have done all of this while largely
preserving the integrity of the underlying 
set-valued data, requiring only small perturbations while
providing strong privacy guarantees to the data of interest, letting
us simultaneously keep data private and preserve its usefulness
in applications.

\section{Conclusion} \label{sec:conclusion}
Future work includes broadening the types of privacy that can be provided
to sets. In this work, we deliberately preserved the shape of sets when
keeping them private, though in other applications it may be desirable
to distort the boundaries of sets, and future work includes defining
privacy of this kind. Similarly, the orientation of a set may be revealing,
and another direction of future work includes protecting a set's orientation
using differential privacy. 
Beyond that, many applications have been characterized
when Gaussian disturbances are added, and developing a Gaussian Perturbation Mechanism
would let users keep set-valued data private while relying on existing convergence
and stability guarantees that already accommodate Gaussian noise.

\bibliographystyle{plain}{}
\bibliography{sources}
\end{document}